# Method of Separating Tangents

In the book [12] on page 218, the following inequality was proven:

$$a^4 + b^4 + c^4 \geq 4abc - 1, \tag{1}$$

where $a$, $b$, $c$ are real numbers. The proof can be written in one line:

$$a^4 + b^4 + c^4 - 4abc + 1 = (a^2 - 1)^2 + (b^2 - c^2)^2 + 2(bc - a)^2 \geq 0.$$

It is interesting to consider, how might one deduce the above on his or her own? Unfortunately, there are no universal ways to identify the above representation. The inevitable question thus arises, can a similar inequality be proved using an alternative, more general method? Such a method of proving the above inequality is presented below.

*Proof.* Let us consider a function: $f(x) = x^4$, where $x \in R$. We first solve the following equation for functions of type $g(x) = kx^3 + m$, at a point such that $f(1) = g(1)$, $f'(1) = g'(1)$ and $f(x) \geq g(x)$, where $x$ is any real number. This demands that numbers $k$ and $m$ be governed by the equation $k + m = 1$, $4 = 3k$. Hence, $k = \frac{4}{3}$, $m = -\frac{1}{3}$. Now it is left to verify that $f(x) \geq g(x)$ for all values of $x \in R$. The latter is true since the inequality $x^4 \geq \frac{4x^3 - 1}{3}$ is (see Pic. 1) equivalent to $(x-1)^2 \left(3\left(x + \frac{1}{3}\right)^2 + \frac{2}{3}\right) \geq 0$, which is evidently true for all values of $x$. Therefore, using Cauchy's inequality $|a|^3 + |b|^3 + |c|^3 \geq 3|abc|$ and $|t| \geq t$ we can write

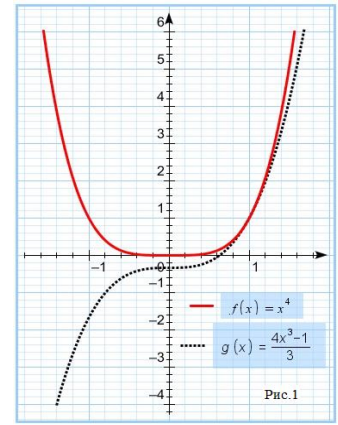

Рис.1

$$a^4 + b^4 + c^4 = |a|^4 + |b|^4 + |c|^4 \geq \frac{4(|a|^3 + |b|^3 + |c|^3)}{3} - 1 \geq 4|abc| - 1 \geq 4abc - 1.$$

Using the above method, the inequality is proven. Is this method better or worse than the original proof? On the one hand, it is much longer. However, this method can be used to solve other inequalities. When doing this, one must analyze the function that is presented in an inequality rather than invent nonstandard transformations. For example, at the Baltic Mathematics Olympiad in 2011, the problem of proving the following inequality was presented:

$$\frac{a}{a^3 + 8} + \frac{b}{b^3 + 8} + \frac{c}{c^3 + 8} + \frac{d}{d^3 + 8} \leq \frac{4}{9}, \tag{2}$$

where $a$, $b$, $c$, $d$ are positive numbers whose aggregate sum is 4.

The solution to (2) was published on the official page of the Baltic Olympiad (see [20]), which can be written in a short form as:

$$\frac{a}{a^3 + 8} + \frac{b}{b^3 + 8} + \frac{c}{c^3 + 8} + \frac{d}{d^3 + 8} =$$

$$= \frac{a}{a^3 + 1 + 1 + 6} + \frac{b}{b^3 + 1 + 1 + 6} + \frac{c}{c^3 + 1 + 1 + 6} + \frac{d}{d^3 + 1 + 1 + 6} \leq [AM - GM] \leq$$

$$\leq \frac{a}{3a+6} + \frac{b}{3b+6} + \frac{c}{3c+6} + \frac{d}{3d+6} = \frac{4 - 2\left(\frac{1}{a+2} + \frac{1}{b+2} + \frac{1}{c+2} + \frac{1}{d+2}\right)}{3} \leq$$

$$\leq [AM - HM] \leq \frac{4 - 2 \cdot \frac{16}{a+2+b+2+c+2+d+2}}{3} = \frac{4}{9}.$$

The following standard abbreviations were used:

AM-GM – Difference between arithmetic mean and geometric mean of numbers;
AM-HM – Difference between arithmetic mean and harmonic mean of numbers.

Once again a question arises, can one can easily arrive at the above solution and can the problem can be solved without using special cases similar to 8=1+1+6? First, it is worth noting that the Baltic Olympiad is a competition between teams with many countries are represented, such as Russia, Germany, Latvia, Finland, Sweden and Estonia. In 2011, of the 11 participating teams, only five were able to solve the above problem, and in the 2012 International Mathematics Olympiad only 189 out of 548 participants were able to solve the same problem. In the author's opinion, such a low success rate has to do with the fact that classical knowledge of inequalities is not sufficient. It is therefore essential for us to develop more general and simpler methods for solving inequalities, at least for certain types commonly encountered. To develop such a method we had to generalize it, to provide theoretical grounding. The theoretical justification is presented in this paper. Based on our results, we were able to develop a new method for solving inequalities using separating tangents. The description of the method and its applications were accepted by the following four publications:

1) Cheboksary, Russia

"About one method of solving inequalities" by I.Z. Ibatulin, A.N. Lepes,
// Works of XXI conference «Mathematics. Education» in 2013, 16 pages (pages. 34–50);

2) Sofia, Bulgaria

"Application of the method of separating tangents to prove inequalities by Ibatulin I.Z., Lepes A.N.
// Didactical Modeling: e-journal 2013. URL: http://www.math.bas.bg/omi/DidMod/index.htm. 13 pages (to be published in 2014);

3) Hong Kong, China

Ibatulin I.Zh., Lepes A.N. Using tangent lines to prove inequalities (part II)
// Mathematical Excalibur. 2013–2014. V.18 N.4, 6 pages.

4) Moscow, Russia

Ibatulin I.Zh., Lepes A.N. Using tangent lines to prove inequalities
// Mathematics in School. 2014. №4, p. 20–24.

The readers are encouraged to solve the following sample problems on their own in order to fully appreciate the practical application of the method of separating tangents.

**Sample Problem 1.** ([16], page. 31, ex. 1.2.9) For any positive numbers $a$, $b$, $c$ such that $a^2+b^2+c^2=3$, prove the inequality: $\frac{1}{a^3+2} + \frac{1}{b^3+2} + \frac{1}{c^3+2} \geq 1$.

**Sample Problem 2.** ([16], page. 55, ex. 3.1.4) For any negative numbers $a$, $b$, $c$, $d$, $e$ such that

prove that
$$\frac{1}{4+a}+\frac{1}{4+b}+\frac{1}{4+c}+\frac{1}{4+d}+\frac{1}{4+e}=1,$$

$$\frac{a}{4+a^2}+\frac{b}{4+b^2}+\frac{c}{4+c^2}+\frac{d}{4+d^2}+\frac{e}{4+e^2}\leq 1.$$

**Sample Problem 3.** (China, 2005, [12], page. 196, problem 132) For any given positive numbers $a$, $b$, $c$, the aggregate sum of which is 1, prove that
$$10(a^3+b^3+c^3)-9(a^5+b^5+c^5)\geq 1.$$

**Sample Problem 4.** ([12], page. 196, problem 124) For any given positive numbers $a$, $b$, $c$ such that $a^2+b^2+c^2=12$, find the largest possible value of the following expression
$$A=a\cdot\sqrt[3]{b^2+c^2}+b\cdot\sqrt[3]{c^2+a^2}+c\cdot\sqrt[3]{a^2+b^2}.$$

**Sample Problem 5.** (Baltic Way, 2002, problem 4, [21]) Let $n$ be a positive whole number. For any numbers of $x_1, x_2, \ldots, x_n \geq 0$ such that $\sum_{k=1}^{n} x_k = 1$, prove the following inequality
$$\sum_{k=1}^{n} x_k(1-x_k)^2 \leq \left(1-\frac{1}{n}\right)^2.$$

### Method of Separating Tangents

We present solutions to the above five sample problems. However, it is worth noting that the application of the method is not always straightforward, and in some cases can be indirect. For example, although Lemma 1 does not provide a formula to identify the local base curve, it does contain conditions for such curve fitting, which can be identified given that the type of the curve is well known. This is similar to Fermat's differential theorem in which, although there are no formulas to identify local extreme points, ways are provided for their identification, and this theorem is widely used in everyday school practice. Using Lemma 2, the inequality proof between a function and its local base curve can be simplified. Without this lemma one has to turn to more involved approaches, such as using Young's inequality, which was the tactic of School number 239 in problem №1 of the Open Olympiad in Saint Petersburg (see [2], page 71).

Here, we present one more demonstration of the method of separating tangents by solving the inequality (2) presented previously.

*Proof.* We note that when $a=b=c=d=1$, then $\frac{a}{a^3+8}+\frac{b}{b^3+8}+\frac{c}{c^3+8}+\frac{d}{d^3+8}=\frac{4}{9}$. Let us introduce $f(x)=\frac{x}{x^3+8}$, where $x\in(0;4)$. Also, we write the following equation tangent to the graph of the function $f$ near $x_0=1$:

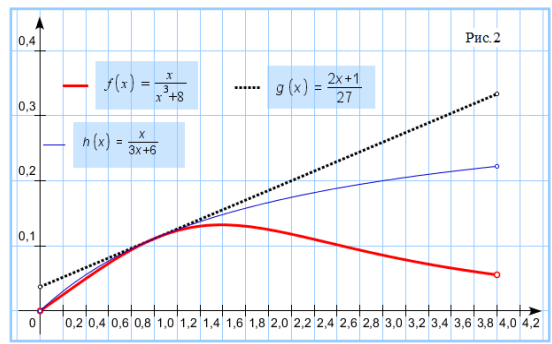

$$y=f(1)+f'(1)(x-1)=\frac{1}{9}+\frac{2}{27}(x-1)=\frac{2x+1}{27}.$$

Let us first show that between (0; 4) the graph of the function $f$ lies below the tangent line of $y = \frac{2x+1}{27}$ (see Pic. 2). Since

$$\frac{x}{x^3 + 8} \leq \frac{2x + 1}{27} \Leftrightarrow (x - 1)^2(2x^2 + 5x + 8) \geq 0, \qquad (3)$$

the above mentioned statement is true. Using the inequality (3), we have

$$\frac{a}{a^3 + 8} + \frac{b}{b^3 + 8} + \frac{c}{c^3 + 8} + \frac{d}{d^3 + 8} \leq \frac{2(a + b + c + d) + 4}{27} = \frac{4}{9}.$$

Thus we have proved inequality (2).

The question we now pose is, how can this tangent line be used, and is there a way to generalize it for several variables?

The school mathematics program defines tangents (see. [11], pages 99–100) as boundary transverse lines. In the article by V.I. Gavrilova (see.[4]), it was pointed out that when a function is downward convex, its graph lies on one side of the tangent line drawn at any point. In this article, the reverse statement was also shown. In the general case, when a function is not necessarily convex, any tangent line has the following interesting property: if the function is continuous, twice differentiable and its second derivative at point $x_0$ is nonzero, then there exists a region near $x_0$ where the graph of the function lies on one side of the tangent line drawn at point $x_0$. In other words, when the second derivative of the function at point $x_0$ is positive then, owing to the continuity of the second derivative, there exists a region near point $x_0$ in which the second derivative is positive, and hence the function is downward convex. The graph of such a function will not lie below the segment of the tangent line near $x_0$.

One of the widely accepted properties of convex functions is that they comply with Jensen's inequality (see, example [5]):

$$\frac{\sum_{k=1}^{n} f(x_k)}{n} \geq f\left(\frac{\sum_{k=1}^{n} x_k}{n}\right),$$

where $f$ is downward convex in the interval $I$, $x_1, x_2, \ldots, x_n \in I$. Problems in mathematics Olympiads contain inequalities related to Jensen's inequalities, which can be written in the following general form:

$$\sum_{j=1}^{n} f(x_j) \geq n f(x_0), \qquad (4)$$

where $x_0, x_1, x_2, \ldots, x_n \in I$, $\sum_{j=1}^{n} l(x_j) = n \cdot l(x_0)$. Herein, when the function $f$ meets the requirement of inequality (4), we will say that the function complies with Jensen's inequality at point $x_0$.

We have observed many problems in mathematical Olympiads related to nonconvex functions that comply with inequality (4). Is this a coincidence, or is there a condition more general than just convexity for inequality (4)? The significance of this question lies the fact that there are at least six different solutions to the one inequality (see, for example, [2], pages 71–75), and a lack of understanding of this may lead to an incorrect deduction of convexity for functions that comply with Jensen's inequality (for example, see [17], p. 19).

Problem №1 of the 2011 Open Olympiad held in School Number 239, Saint Petersburg (see [2], page 71) is defined as:

$$\frac{a}{a^3+4}+\frac{b}{b^3+4}+\frac{c}{c^3+4}+\frac{d}{d^3+4}\leq\frac{4}{5}, \qquad (5)$$

where *a*, *b*, *c*, *d* are positive numbers whose aggregate sum is 4. A key point in the solution to this problem is the following:

$$\frac{a}{a^3+4}\leq\frac{2a+3}{25}.$$

The above inequality may only be found using Young's inequality. However, it turns out that the line $y=\frac{2x+3}{25}$ is a tangent to the graph of the function $f(x)=\frac{x}{x^3+4}$, where $x\in(0; 4)$, and the difference $\frac{a}{a^3+4}-\frac{2a+3}{25}$ in the range of (0; 4) has the same polarity as the polynomial $(a-1)^2(-2a^2-5a-8)$. Therefore, the supporting inequality for the solution of problem №1 could have been proven by using factorization.

In the solution to problem №4 ([2], pages 73–74), called "Jensen's inequality enhanced by a stick", because the function *f* is non-convex in the entire range of (0; 4), Jensen's inequality is only applied over some part of this range, and for the remaining part verbal proof is used. Moreover, it turns out that even though the function is non-convex, the Jensen's inequality is true at $x_0$=1. Thus the enhancement of Jensen's inequality could have been avoided because the graph of the function *f* does not lie below the tangent line drawn to it at point $x_0$=1.

Is it a coincidence that the point $x_0$=1 is the double root of the polynomial numerator of a fraction, that is similar to $\frac{a}{a^3+4}-\frac{2a+3}{25}$? Or is it a property of difference between a function and its tangent line at the point $x_0$? These and other questions are addressed in this paper.

With regard to the possible incorrect deduction of convexity (see, [17], page 19), even though the function $f(x)=10x^3-9x^5$ between [0; 1] is non-convex (see Pic. 5), it still complies with Jensen's inequality at the point $x_0$=1/3:

$$10(a^3+b^3+c^3)-9(a^5+b^5+c^5)=f(a)+f(b)+f(c)\geq 3f\left(\frac{1}{3}\right)=1,$$

where *a*, *b*, *c*>0 and *a*+*b*+*c*=1. The latter inequality is shown using tangents in Theorem 2 below and without using tangents in [12], page 297.

## Theory

In most literature dedicated to Jensen's inequality only convex functions are mentioned (see, e.g. [3], [5], [6]). However, there are also expectations (see [16], pages 68–70). This paper attempts to fill this gap. We suggest expanding our definition of tangents in the following way:

*Definition* 1. We introduce a numerical function *f*, defined in some interval *I*. Let continuous function *g*: *I*→R be local base to the graph of a function *f* at a point $x_0 \in I$. If $f(x_0)=g(x_0)$ then there exists $\delta$>0 such that one of the following conditions are true:

1) For any $x \in I \cap (x_0-\delta; x_0+\delta)$, inequality $f(x) \geq g(x)$ is true;

2) For any $x \in I \cap (x_0-\delta; x_0+\delta)$ inequality $f(x) \leq g(x)$ is true.

Observe that if the function is continuous and twice differentiable and its second derivative at the point $x_0$ is not 0, then the tangent line complies with Condition 1. We cannot say in general though that tangent lines comply with the Condition 1. For instance, the tangent line to the function $f_2(x) = x^3$ at the point $x_0=0$ does not comply.

**Lemma 1 (special case of Fermat's Theorem).** Let us introduce numerical functions $f$ and $g$, defined in the interval $I$. Assume $f$ and $g$ are differentiable at point $x_o \in I$ and that $f(x_0) = g(x_0)$. If there exist $\delta > 0$ such that, for any $x \in I \cap (x_0 - \delta; x_0 + \delta)$ inequality $f(x) \geq g(x)$ is true, then $f'(x_0) = g'(x_0)$.

*Proof.* Based on the theorem we can write

$$h(x) = f(x) - g(x) \geq 0 = h(x_0), \text{ where } x \in I \cap (x_0 - \delta; x_0 + \delta).$$

Here, the point $x_0$ is a local minimum point of the function $h$. Thus based on Fermat's theorem, since $h'(x_0) = 0$, then $f'(x_0) = g'(x_0)$. The proof is complete.

The practical point of Lemma 1 is that for a differentiable function $f$ at a point $x_0$, the local base curve must be searched from those functions whose derivative at $x_0$ is the same as $f'(x_0)$. Therefore, Lemma 1 explains the recommendations provided by Hung in his book (see [16], page 136). In certain cases, this can be quite useful.

When proving inequality (1), and also inequality (2), we studied the difference between the function $f(x) = x^4$, where $x \in R$, and its local base curve $g(x) = \frac{4x^3-1}{3}$; the following polynomial was identified

$$(x-1)^2 \left( 3\left(x+\frac{1}{3}\right)^2 + \frac{2}{3} \right),$$

whose intersection point coincided with its twice repeating root.

**Lemma 2.** *Let us introduce polynomial functions P, Q and g, in the interval I, where $f(x_0) = g(x_0)$ and $f'(x_0) = g'(x_0)$ at the point $x_0 \in I$, and $f(x) = \frac{P(x)}{Q(x)}$. Then the difference $h(x) = f(x) - g(x)$ can be written as $\frac{(x-x_0)^2 \cdot T(x)}{Q(x)}$, where T is a polynomial.*

*Proof.* Based on Taylor's expansion, any polynomial can be written in the following way:

$$P(x) = \sum_{k=0}^{n} \frac{P^{(k)}(x_0)}{k!}(x - x_0),$$

where $n$ is the degree of a polynomial $P$, and $P^{(k)}(x_0)$ is the $k$-th derivative of polynomial $P$ at the point $x_0$. Hence, polynomials $P, Q, g$ can be written as

$$P(x) = P(x_0) + P'(x_0)(x - x_0) + P_1(x)(x - x_0)^2,$$

$$Q(x) = Q(x_0) + Q'(x_0)(x - x_0) + Q_1(x)(x - x_0)^2,$$

$$g(x) = \frac{P(x_0)}{Q(x_0)} + g_1(x)(x - x_0),$$

where $P_1$, $Q_1$, $g_1$ – are polynomials. Therefore,

$$f(x) - g(x) = \frac{P(x)}{Q(x)} - \frac{P(x_0)}{Q(x_0)} - g_1(x)(x - x_0) =$$

$$= \frac{(x-x_0)\left(P'(x_0)Q(x_0)-P(x_0)Q'(x_0)-g_1(x)Q(x)Q(x_0)+(x-x_0)(P_1(x)Q(x_0)-P(x_0)Q_1(x_0))\right)}{Q(x)Q(x_0)}.$$

Now, $f'(x_0) = g'(x_0)$ implies that $g_1(x) = f'(x_0) + g_2(x)(x - x_0)$, where $g_2$ is a polynomial. Hence,

$$f(x) - g(x) = \frac{(x-x_0)^2\left[(P_1(x)Q(x_0)-P(x_0)Q_1(x))-f'(x_0)Q(x_0)\left(Q'(x_0)+Q_1(x)(x-x_0)\right)-g_2(x)Q(x)Q(x_0)\right]}{Q(x)Q^2(x_0)}.$$

The proof is complete.

*Solution of Sample Problem 1.* We note that when $a = b = c = 1$ the provided inequality becomes an equation. Let $f(x) = \frac{1}{x^3+2}$, $g(x) = kx^2 + m$, where $x \in (0; \sqrt{3})$. Numbers $k$ and $m$ are chosen such that $f(1) = g(1)$, $f'(1) = g'(1)$. Hence, $k + m = \frac{1}{3}$, $2k = -\frac{1}{3}$. Therefore, $g(x) = -\frac{x^2}{6} + \frac{1}{2}$. Inequality $\frac{1}{x^3+2} \geq -\frac{x^2}{6} + \frac{1}{2}$ is true (see Pic. 3), since it is equivalent to the inequality $x^2(x - 1)^2(x + 2) \geq 0$. Hence,

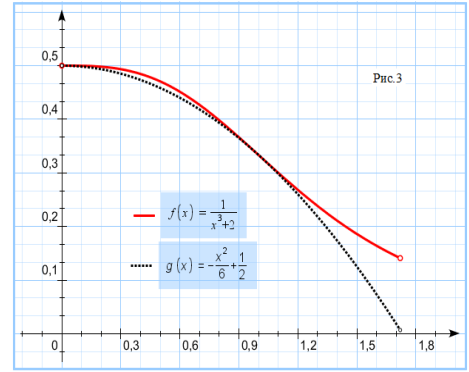

Рис. 3

$$\frac{1}{a^3+2} + \frac{1}{b^3+2} + \frac{1}{c^3+2} \geq \frac{3}{2} - \frac{a^2+b^2+c^2}{6} = 1.$$

Thus, the proof is complete.

It is interesting to note that the graph of the function *f* not being below its local base curve is a sufficient condition for the inequality (4).

**Theorem 1 (sufficient conditions for Jensen's inequality).** *Let f and l be numerical functions defined in the interval I. Assume that functions f and l are differentiable at the point $x_o \in I$. If for every $x \in I$ the following inequality is true:*

$$f(x) \geq k \cdot l(x) + m,$$

*where* $= \begin{cases} 0, \text{ when } l'(x_0) = 0, \\ \frac{f'(x_0)}{l'(x_0)}, \text{when } l'(x_0) \neq 0, \end{cases}$ $m = f(x_0) - k \cdot l(x_0),$

*then Jensen's inequality is true for the function f at the point $x_0$.*

*Proof* of Theorem 1 is trivial from the condition of the theorem

$$\sum_{j=1}^{n} f(x_j) \geq \sum_{j=1}^{n} (k \cdot l(x_j) + m) = k \sum_{j=1}^{n} l(x_j) + mn = n(k \cdot l(x_0) + m) = nf(x_0)$$

*Note 1.* If $l(x) = x$ and the function $f$ is downward convex on the interval $I$, then the graph of the function lies on one side of the tangent line drawn at any point $x_0 \in I$, so that the following inequality $f(x) \geq kx + m$ is true where $k$ and $m$ are defined as in Theorem 1. Consequently, based on Theorem 1, for any convex function Jensen's inequality is true. Jensen's inequality is true at point $x_0$ for any convex or non-convex function, the graph of which lies on one side of the tangent line drawn at the point $x_0$. This is a general statement that is mentioned in, for example, article [5].

*Solution of Sample Problem 2.* Let $f(x) = \frac{x}{4+x^2}$, $g(x) = \frac{k}{4+x} + m$, where $x \geq 0$ and $k$ and $m$ are such that $f(1) = g(1)$, $f'(1) = g'(1)$. Hence, $k=-3$, $m=0{,}8$. Since inequality $\frac{x}{4+x^2} \leq \frac{4}{5} - \frac{3}{4+x}$ is equivalent to inequality $(x-1)^2(x+1) \geq 0$, then it is true for any $x \geq 0$ (see Pic. 4). Therefore,

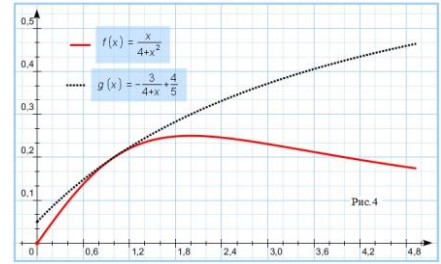

$$\frac{a}{4+a^2} + \frac{b}{4+b^2} + \frac{c}{4+c^2} + \frac{d}{4+d^2} + \frac{e}{4+e^2} \leq 4 - 3\left(\frac{1}{4+a} + \frac{1}{4+b} + \frac{1}{4+c} + \frac{1}{4+d} + \frac{1}{4+e}\right) = 1.$$

The proof is complete.

Theorem 1 can be generalized in the following way.

**Theorem 2 (sufficient conditions for Jensen's inequality).** *Let $f$ and $l$ be numerical functions defined in the interval I, and G is a set, $x_0 \in I/G$. Assume functions $f$ and $l$ are differentiable at the point $x_0$ and at the set G and I the function reaches its minimum value. If*

$$\min_G f + (n-1) \min_I f \geq nf(x_0)$$

*and for any $x \in I/G$ the following inequality is true,*

$$f(x) \geq k \cdot l(x) + m,$$

*where $= \begin{cases} 0, \text{when } l'(x_0) = 0, \\ \frac{f'(x_0)}{l'(x_0)}, \text{when } l'(x_0) \neq 0, \end{cases}$ $m = f(x_0) - k \cdot l(x_0)$,*

*then Jensen's inequality is true for the function $f$ at the point $x_0$.*

*Proof.* If $x_1, x_2, \ldots, x_n \in I/G$, then, going through similar steps as in Theorem 1, we find that $\sum_{j=1}^{n} f(x_j) \geq nf(x_0)$.

Let us allow, for example, that $x_1$ does not belong to the set $I/G$, i.e. $x_1 \in G$. Then according to the definition of the minimum value of the function we have

$$f(x_1) \geq \min_G f, f(x_2) \geq \min_I f, f(x_3) \geq \min_I f, \ldots, f(x_n) \geq \min_I f.$$

Adding the obtained inequalities we arrive at the conditions set by the theorem.

Since the process of selecting a set $G$ and identifying the minimum value of a function requires some work, after which arriving at the results is straightforward, the proof of the Theorem 2 is more useful than its result. Thus, by using Theorem 2 and the graph of a function it can be decided if a local base curve can be used or not. However, the full proof of the presented inequalities requires similar steps to those taken in Theorem 2.

*Solution of Sample Problem* 3. Let $f(x) = 10x^3 - 9x^5$, $x \in (0; 1]$. When $a = b = c = \frac{1}{3}$ we obtain the following equation

$$10(a^3 + b^3 + c^3) - 9(a^5 + b^5 + c^5) = 1.$$

Let us plot a tangent line to the graph of the function $f$ at the point $x_0 = \frac{1}{3}$:

$$y = f\left(\frac{1}{3}\right) + f'\left(\frac{1}{3}\right)\left(x - \frac{1}{3}\right) = \frac{1}{3} + \frac{25}{9}\left(x - \frac{1}{3}\right) = \frac{25}{9}x - \frac{16}{27}.$$

The following inequality (see Pic. 5)

$$10x^3 - 9x^5 \geq \frac{25}{9}x - \frac{16}{27}. \qquad (6)$$

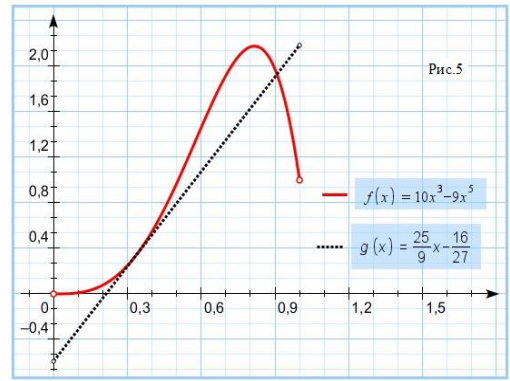

can be also rewritten as

$$\left(x - \frac{1}{3}\right)^2 \left(9x^3 + 6x^2 - 7x - \frac{16}{3}\right) \leq 0.$$

Since for any $x \in (0; 0{,}9)$ it is true that

$$9x^3 + 6x^2 - 7x - \frac{16}{3} = (x - 0{,}9)(9x^2 + 14{,}1x + 5{,}69) - \frac{637}{3000} < 0,$$

inequality (6) is true for all $x \in (0; 0{,}9)$. Therefore, if positive numbers $a$, $b$, $c$ are less than 0,9 and their aggregate sum is 1, then from inequality (6) we can write

$$10(a^3 + b^3 + c^3) - 9(a^5 + b^5 + c^5) \geq \frac{25}{9}(a + b + c) - 3 \cdot \frac{16}{27} = 1.$$

However, if any of $a$, $b$, $c$ are not less than 0,9, then because function $f$ is decreasing in the interval [0,9; 1] since

$$f'(x) = -45x^2\left(x^2 - \frac{2}{3}\right) < 0, where\ x \in [0{,}9; 1],$$

and because function $f$ is nonnegative on the entire interval of $(0; 1]$, we can write

$$10(a^3 + b^3 + c^3) - 9(a^5 + b^5 + c^5) = f(a) + f(b) + f(c) \geq$$

$$\geq \min_{[0{,}9;\ 1]} f(x) = f(1) = 1.$$

Therefore, the proof is complete.

One of the drawbacks of Theorems 1 and 2 is that their type of base curves and conditions for variables are the same. To avoid this, we could use properties of power means (for more details see [7], pages 12–23).

*Definition* 2. Let us define positive $x_1$, $x_2$, …., $x_n$. When $\alpha$ is not zero, the power mean of $x_1$, $x_2$, …., $x_n$ is defined as

$$c_\alpha(x_1, x_2, \ldots, x_n) = \left(\frac{\sum_{j=1}^{n} x_j^\alpha}{n}\right)^{\frac{1}{\alpha}}.$$

When $\alpha=0$,

$$c_0(x_1, x_2, \ldots, x_n) = \sqrt[n]{x_1 \cdot x_2 \cdot \ldots \cdot x_n}.$$

**Theorem 3 (sufficient conditions for Jensen's inequality).** *For $\alpha$, positive number $x_0$ and whole number $n \geq 2$, let a function f be defined on all sets of positive numbers and be differentiable at the point $x_0$. If $(\alpha - 1) \cdot f'(x_0) \leq 0$ and for any positive x such that $x^\alpha < nx_0^\alpha$, and the following inequality $f(x) \geq f(x_0) + f'(x_0)(x - x_0)$ is true, then*

$$\sum_{j=1}^{n} f(x_j) \geq nf(x_0),$$

where $x_1$, $x_2$, …, $x_n > 0$ and $c_\alpha(x_1, x_2, \ldots, x_n) = x_0$.

*Proof.* When $\alpha \neq 0$ and for positive $x_1$, $x_2$, …, $x_n > 0$, $c_\alpha(x_1, x_2, \ldots, x_n) = x_0$ is true, then regardless of whether the sign of $\alpha$ is positive or negative, for any $k=1, 2, \ldots, n$ the inequality $x_k^\alpha < nx_0^\alpha$ is true. When $\alpha=0$, the inequality is also true.

We assume $\alpha<1$. Then $f'(x_0) \geq 0$. Since power means are monotonic (see example [7], page 21), we can write

$$\sum_{j=1}^{n} f(x_j) \geq \sum_{j=1}^{n} \left(f(x_0) + f'(x_0)(x_j - x_0)\right) = nf(x_0) + nf'(x_0)(c_1(x_1, x_2, \ldots, x_n) - x_0) \geq$$

$$\geq nf(x_0) + nf'(x_0)(c_\alpha(x_1, x_2, \ldots, x_n) - x_0) = nf(x_0)$$

The case when $\alpha>1$ is equivalent to the above, and when $\alpha=1$ it follows from Theorem 1. Therefore, the proof of Theorem 3 is complete.

*Note* 2. When $\alpha=1$ in Theorem 3, it follows that for any downward convex function the Jensen's inequality is true at a point $x_0$.

Theorem 3 has important practical applications. When the aggregate sum of squares of variables is constrained, then the local base curve can be applied if the derivative of a function at a point is negative. When the product of variables is fixed, then the derivative of a function at a point should be positive. However, if the conditions for the function are not met, then a polynomial having the same degree as a function can be selected as a tangent. Therefore, when the aggregate sum of squares of variables is fixed and derivative of a function is positive, then the parabola at a point can be selected as a tangent.

*Solution of Sample Problem* 4. Let $f(x) = -x \cdot \sqrt[3]{12 - x^2}$, where $x \in (0; 2\sqrt{3})$. Let us construct equations for a tangent to the graph $f$ at the point $x_0=2$:

$$y = f(x_0) + f'(x_0)(x - x_0) = -4 - \frac{2}{3}(x - 2) = -\frac{4x+4}{3}.$$

If for any positive $x$ the inequality $f''(x) = \frac{\frac{2}{3}x}{\sqrt[3]{(12-x^2)^2}} + \frac{\frac{16}{9}x^3}{\sqrt[3]{(12-x^2)^5}} > 0$ is true, then the function $f$ is downward convex on the interval $(0; 2\sqrt{3})$, so that its graph does not lie below its tangent (see Pic. 6). Hence, according to Theorem 4 when $\alpha = 2, x_0 = 2, f(2) = -4, f'(2) = -2/3$ we can write

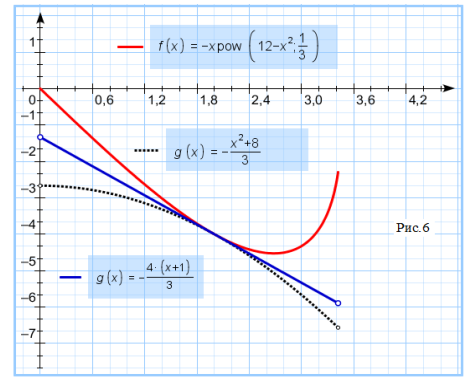

Рис.6

$$-a \cdot \sqrt[3]{b^2 + c^2} - b \cdot \sqrt[3]{c^2 + a^2} - c \cdot \sqrt[3]{a^2 + b^2} \geq -3 \cdot 4 =$$
$$= 12,$$

and equality of both sides is reached when $a=b=c=2$. Therefore,

$$\max_{\substack{a^2+b^2+c^2=3,\\a,b,c>0}} \left(a \cdot \sqrt[3]{b^2 + c^2} + b \cdot \sqrt[3]{c^2 + a^2} + c \cdot \sqrt[3]{a^2 + b^2}\right) = 12.$$

*Note* 3. When solving sample problem 4 we used a tangent, but we could have also used a parabola (see Pic. 6). When solving Problem 1 though, we did not consider a tangent line, even though the derivative at the intersecting point was negative, which meets the condition of Theorem 3. It can be deduced from Picture 3 that the graph of the function does not lie on one side of the tangent line, so that the domain of the function had to be partitioned as it was done when solving Problem 3. It is evident that even though such a proof was feasible, it would be more complicated than the one presented above.

Theorem 3 can be used to identify when a tangent line can be utilized, but also when other curves should be selected. For example, in Problem 5 of the 2005 Baltic Olympiad, the following inequality was presented:

$$\frac{a}{a^2 + 2} + \frac{b}{b^2 + 2} + \frac{c}{c^2 + 2} \leq 1,$$

where *a*, *b*, *c* are positive such that *abc*=1.

From the 11 participating teams, only Estonia, Iceland and Poland could successfully solve the above problem. The book [12] (page 246, problem 56) contains a solution where Cauchy's inequality is used. To effectively use Cauchy's inequalities the representations of type $a^2+2=a^2+1+1$ must be used, which is trivial, but unimaginable as to how one might easily arrive at them. In fact eight of the 11 teams could not arrive at this representation. Let us look at this problem using Theorem 3. The derivative of the function $f(x) = -\frac{x}{x^2+2}$, where $x>0$, at the point $x_0=1$, must be positive since α=0. But, in this case the derivative is negative. Therefore, consideration of a tangent line as a base curve is not required, even though the function *f* is downward convex; other types of curves should be used. Moreover, based on similar cases solved by Hung [16] (page 136), we could use the following function

$$g(x) = -\frac{3 + \ln x}{9}, \text{ where } x > 0.$$

The graph of the function *f* does not lie on one side of the graph of the function *g*, but by undertaking similar steps as in Problem 3, the presented inequality can be solved.

There are, however, some complications in Theorem 3 with regard to solving power function inequalities by using tangent lines as local base curves. Furthermore, power functions could be also used as base curves of functions with constrained aggregate sum of variables.

**Theorem 4 (sufficient conditions for Jensen's inequality).** *Let $\alpha \neq 0$, $x_0$ be a positive number, and $n \geq 2$ a whole number. Let the function f be defined for all sets of positive numbers and differentiable at the point $x_0$. If $(\alpha - 1) \cdot f'(x_0) \geq 0$ and for any positive number x such that $x < nx_0$, the inequality $f(x) \geq \frac{f'(x_0)}{\alpha \cdot x_0^{\alpha-1}}(x^\alpha - x_0^\alpha) + f(x_0)$ is true, then*

$$\sum_{j=1}^{n} f(x_j) \geq nf(x_0),$$

*where $x_1, x_2, ..., x_n > 0$ и $x_1 + x_2 + ... + x_n = nx_0$.*
 Proof of Theorem 4 is similar to the proof of Theorem 3.
 Both Theorem 3 and 4 can be generalized as was done for Theorem 1.

**Special Case**

**Theorem 5.** *Let $P(x) = ax^3 + bx^2 + cx + d$ be polynomials, $a \neq 0$, n be a whole number and $x_0$ a positive number. We introduce non-negative numbers $x_1, x_2, ..., x_n$, whose aggregate sum is $nx_0$. If $2ax_0 + b \geq 0$ and $(n+2)ax_0 + b \geq 0$, then Jensen's inequality is true at the point $x_0$.*
 *Proof.* According to Lemma 2, the function $P(x) - P(x_0) - P'(x_0)(x - x_0)$ can be rewritten as $(x - x_0)^2 P_1(x) = (x - x_0)^2 (ax + 2ax_0 + b)$. Since

$$\sum_{j=1}^{n} P'(x_0)(x_j - x_0) = P'(x_0)\left(\sum_{j=1}^{n} x_j - nx_0\right) = 0,$$

the inequality $\sum_{j=1}^{n} P(x_j) \geq nP(x_0)$ is equivalent to the following inequality:

$$\sum_{j=1}^{n}(x_j - x_0)^2 (ax_j + 2ax_0 + b) \geq 0.$$

The aggregate sum of non-negative numbers $x_1, x_2, \ldots, x_n$ is $nx_0$, which implies that each of these numbers lies within the interval of $[0; nx_0]$. However, as the function $ax + 2ax_0 + b$ is monotonic and since $2ax_0+b$ and $(n+2)ax_0+b$ are non-negative, for every $x \in [0; nx_0]$ the inequality $ax + 2ax_0 + b \geq 0$ is true. The latter also implies that the above inequality is also true. Hence, the proof is complete.

*Note 4.* In fact, for the polynomial $P(x)=ax^3+bx^2+cx+d$ to be downward convex, the required condition is $P''(x) = 6ax + 2b \geq 0$ on the entire interval $[0; nx_0]$. Henceforth, the required conditions are $b \geq 0$ and $3nax_0+b \geq 0$. If these conditions are true for $a$, $b$, $n$, $x_0$, then the inequalities $2ax_0+b \geq 0$ and $(n+2)ax_0+b \geq 0$ must also be true. These results are in agreement with the widely accepted fact that Jensen's inequalities are true for convex functions. However, when $a$, $b$, $n$, $x_0$ are such that $2ax_0+b \geq 0$ and $(n+2)ax_0+b \geq 0$, it is not necessarily true that $b>0$ and $3nax_0+b>0$. For instance, when $n=3$, $a=1$, $b=-1$, $x_0=1$, then any function of the type $P(x) = x^3 - x^2 + cx + d$ will be nonconvex on the interval $[0; 3]$. However, Jensen's inequality still holds for the function at the point $x_0=1$. Obviously, another set of $n$, $a$, $b$, $x_0$ can be also chosen, for example $a=1$, $b=-x_0$, $x_0>0$, $n \geq 2$. Therefore, we have shown that there can be infinite nonconvex functions for which Jensen's inequality is true at a given point.

*Solution of Sample Problem* 5. When $n \geq 2$, the inequality presented in the problem can be solved using Theorem 5 and inputs $a=-1$, $b=2$, $c=-1$, $d=0$. When $n=1$, the presented inequality takes the form of $0 \leq 0$, which is evidently true.

## Description of the Method of Separating Tangents

The main feature of our method is in identifying appropriate types of tangents, their key properties and intersecting points, which are used to describe local base curves for proving inequalities.

The method described below is used to prove inequalities of the following type:
$$\sum_{k=1}^{n} f_k(x_k) \geq A$$
and inequalities that can be rearranged into this form.

It was shown before, in Lemmas 1 and 2, that tangents should be chosen such that their values and derivatives match at the point of intersection.

Here, we review three cases:

**1st Case.** $f_1 = f_2 = \ldots = f_n = f$ and the graph of function *f* lies not below local base curve.

a) If $x_1+x_2+\ldots+x_n=B$, $x_1, x_2, \ldots, x_n>0$, then a straight line can be used as the local base curve. However, if the base curve exceeds the function *f* at any point, then a power function can be selected as the base curve. Selection of the power function is done according to Theorem 4. In both cases, the intersection point is chosen as $\frac{B}{n}$.

b) If $\alpha \neq 0$ and $x_1^\alpha + x_2^\alpha + \ldots + x_n^\alpha = B$, $x_1, x_2, \ldots, x_n>0$, then local base curves are represented by functions of type $g(x) = kx^\alpha + m$, where coefficients $k$ and $m$ are found from the equations $g(x_0)=f(x_0)$, $g'(x_0)=f'(x_0)$. The intersection point $x_0$ is chosen to be $\left(\frac{B}{n}\right)^{1/\alpha}$. If any part of the base curve exceeds the function *f*, then the local base curve can be represented as a straight line according to Theorem 3.

c) If $x_1 \cdot x_2 \cdot \ldots \cdot x_n = B$, $x_1, x_2, \ldots, x_n > 0$, then local base curves are represented by functions of type $k \cdot \ln x + m$, where $k = x_0 \cdot f'(x_0)$, $m = f(x_0) - k \ln x_0$, $x_0 = \sqrt[n]{B}$. As in 1b), if the local base curve does not comply with the set constraints, then straight lines can be used given that the derivative of the function $f$ at the point $x_0$ is positive.

d) If a function $F(x_1, x_2, \ldots, x_n)$ is homogenous, then for the following inequalities
$$F(x_1, x_2, \ldots, x_n) \geq 0,$$
where $x_1, x_2, \ldots, x_n > 0$, it is sufficient to prove $c_\alpha(x_1, x_2, \ldots, x_n) = 1$. The selection of $\alpha$ can be made according to Theorem 3. The type of local base curve and intersection point are selected as in 1b).

e) If the function $F(x_1, x_2, \ldots, x_n)$ is homogenous of **degree 0**, then for the following inequalities
$$F(x_1, x_2, \ldots, x_n) \geq A,$$
where $x_1, x_2, \ldots, x_n > 0$, $A \neq 0$, it is sufficient to prove $c_\alpha(x_1, x_2, \ldots, x_n) = 1$. The selection of $\alpha$ can be made according to Theorem 3. The type of local base curve and intersection point are selected as in 1b).

f) When there are no constraints on variables or these constraints are different than those shown above, then to prove inequalities of type $F(x_1, x_2, \ldots, x_n) \geq 0$ for any positive $S$, positive numbers $x_1, x_2, \ldots, x_n$ are chosen whose aggregate sum is $S$. Then, the proof continues as in 1a), where constraints on $S$ are found based on generic inequalities.

**2$^{nd}$ Case.** If for functions $f_1, f_2, \ldots, f_n$, it is necessary to prove the inequality $\sum_{k=1}^{n} f_k(x_k) \geq A$, given the constraints on variables $\sum_{j=1}^{n} l(x_j) = B$, then the local base curve can be represented by a function of type $k_j l(x) + m_j$, $k_j = \dfrac{f_j'(x_j^{(0)})}{l'(x_j^{(0)})}$, $m_j = f_j(x_j^{(0)}) - kl(x_j^{(0)})$, $j=1, 2, \ldots, n$, and intersection points $x_1^{(0)}, x_2^{(0)}, \ldots, x_n^{(0)}$ are found from the following system of equations:

$$\sum_{j=1}^{n} l(x_j^{(0)}) = B, \quad \frac{f_1'(x_1^{(0)})}{l'(x_1^{(0)})} = \frac{f_2'(x_2^{(0)})}{l'(x_2^{(0)})} = \ldots = \frac{f_n'(x_n^{(0)})}{l'(x_n^{(0)})}. \tag{6}$$

Here it is sufficient for us to obtain one solution of equation (6). In a more general case $l(x)=x$, the system of equations (6) may be written as:

$$\sum_{j=1}^{n} x_j^{(0)} = B, \quad f_1'(x_1^{(0)}) = f_2'(x_2^{(0)}) = \ldots = f_n'(x_n^{(0)}).$$

**3$^{rd}$ Case.** $f_1 = f_2 = \ldots = f_n = f$ and the graph of function $f$ lies on both sides of the tangent line. In this case, the domain of the function is split into two regions so that in one region the function lies below the tangent. The proof is then continued as in Theorem 2.

For a better representation of cases 1a), 1b) and 1c), the following schematic diagram is drawn for identifying an appropriate local base curve (LBC):

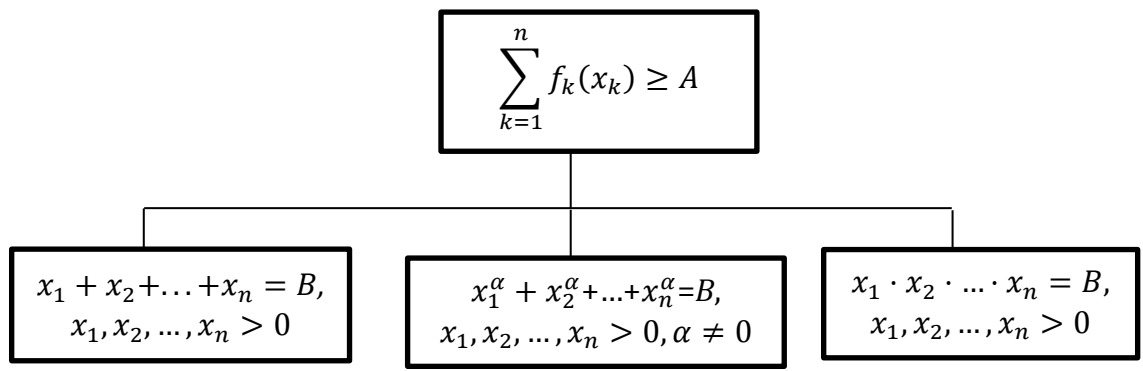

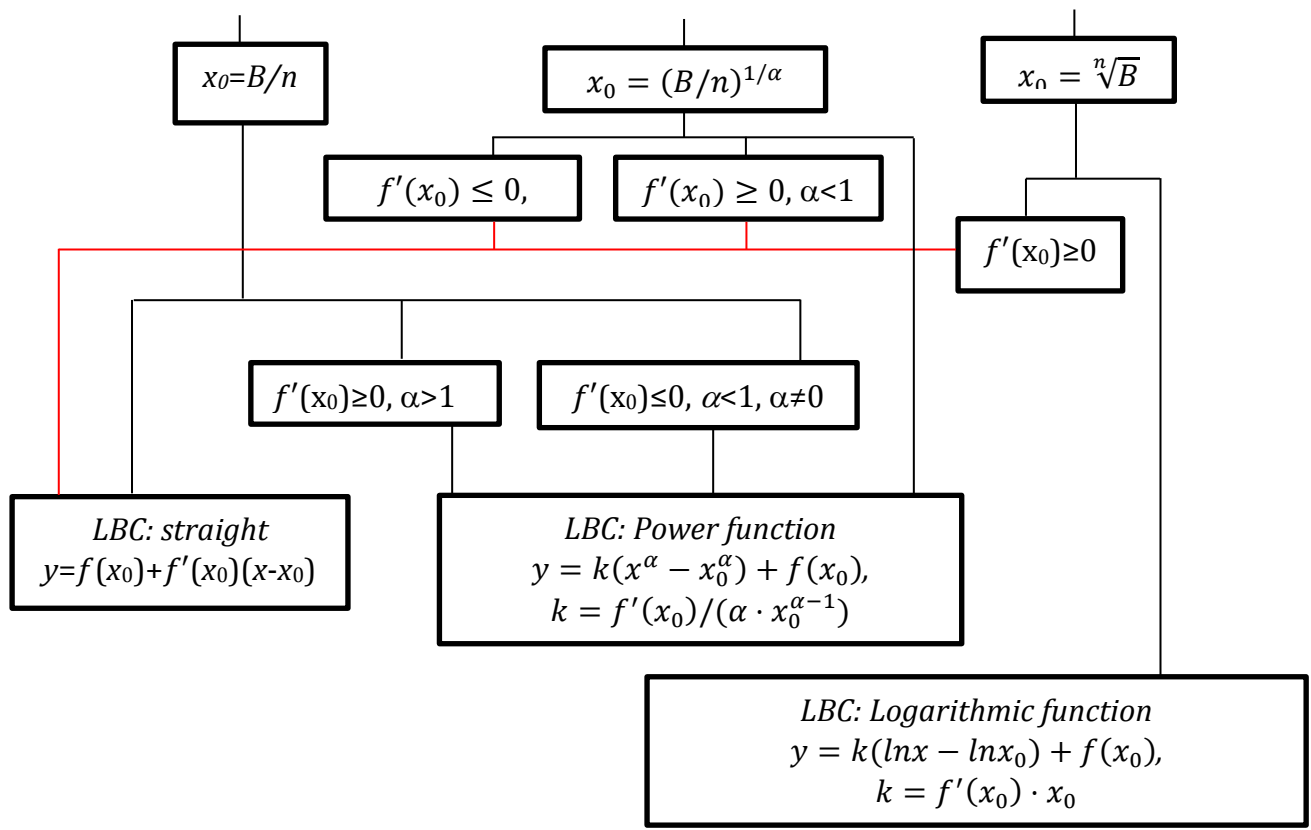

To conclude, we give examples of how the method of separating tangents is applied on several well-known inequality problems.

Book [12], page 189, Problem 63:

**Example 1:** Let *a*, *b*, *c*, *d* be positive numbers such that $a^2+b^2+c^2+d^2=1$. Prove the following inequality:

$$\sqrt{1-a} + \sqrt{1-b} + \sqrt{1-c} + \sqrt{1-d} \geq \sqrt{a} + \sqrt{b} + \sqrt{c} + \sqrt{d}.$$

In this problem on pages 250–251, the author states that for any number $x \in (0; 1)$, the following inequality is true

$$\sqrt{1-x} - \sqrt{x} \geq -\sqrt{2}\left(x - \frac{1}{2}\right).$$

However, it can be easily shown that the latter does not always hold, for example, when $x=0.64$ (see Pic. 7). This incorrect result occurred because the author used the fraction $\frac{1-2x}{\sqrt{1-x}+\sqrt{x}}$ in the proof, which can only be used when the numerator $1-2x$ is positive, and this does not always hold on the interval (0; 1). Hence, the incorrect deduction was unavoidable.

The derivative of the function $f(x) = \sqrt{1-x} - \sqrt{x}$ at the intersection point $x_0=0.5$ complies with the Theorem 3, but the graph of the function lies on both sides of the tangent. Therefore, it would have been impossible to use a tangent line in this case because only the condition in Theorem 2 can be met, as demonstrated in picture 7. However, this only shows that the type of local base curve has been poorly chosen. If, for instance, a parabola were chosen as the local base curve, then the presented inequality could be proven.

*Solution.* Let $f(x) = \sqrt{1-x} - \sqrt{x}$, $g(x) = kx^2 + m$, where $0<x<1$. Numbers $k$ and $m$ are chosen such that $f\left(\frac{1}{2}\right) = g\left(\frac{1}{2}\right)$, $f'\left(\frac{1}{2}\right) = g'\left(\frac{1}{2}\right)$. Hence, for $k$ and $m$ the equations $0 = \frac{k}{4} + m$, $-\sqrt{2} = k$ are true. Therefore, $g(x) = -\sqrt{2}\left(x^2 - \frac{1}{4}\right)$. Now let us show that for any $x \in (0; 1)$ the following inequality is true (see Pic. 7)

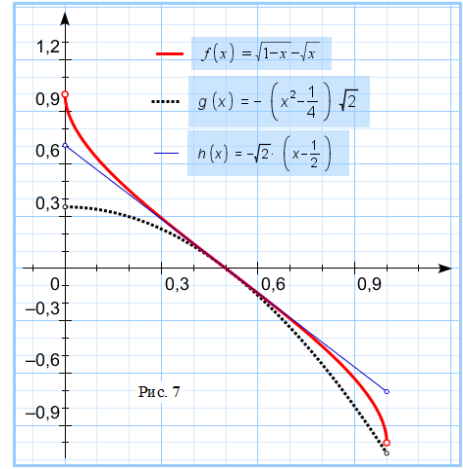

Рис. 7

$$\sqrt{1-x} - \sqrt{x} \geq -\sqrt{2}\left(x^2 - \frac{1}{4}\right) \quad (7)$$

After some rearrangements of the inequality (7), we obtain:

$$\frac{1-2x}{\sqrt{1-x}+\sqrt{x}} \geq \frac{(1-2x)(1+2x)}{2\sqrt{2}} \Leftrightarrow \frac{(1-2x)\left(2\sqrt{2} - (1+2x)(\sqrt{1-x}+\sqrt{x})\right)}{2\sqrt{2}(\sqrt{1-x}+\sqrt{x})} \geq 0 \Leftrightarrow$$

$$\Leftrightarrow \frac{(1-2x)\left(2\sqrt{2} - (1+2x)\left(\frac{1-2x}{\sqrt{1-x}+\sqrt{x}} + 2\sqrt{x}\right)\right)}{2\sqrt{2}(\sqrt{1-x}+\sqrt{x})} \geq 0 \Leftrightarrow$$

$$\Leftrightarrow \frac{(1-2x)\left(\frac{2x-1}{\sqrt{1-x}+\sqrt{x}}(1+2x) + 2(\sqrt{2} - \sqrt{x} - 2x\sqrt{x})\right)}{2\sqrt{2}(\sqrt{1-x}+\sqrt{x})} \geq 0 \Leftrightarrow$$

$$\Leftrightarrow \frac{(1-2x)\left(\frac{2x-1}{\sqrt{1-x}+\sqrt{x}}(1+2x) + 2(1-\sqrt{2x})(\sqrt{2x}+\sqrt{x}+\sqrt{2})\right)}{2\sqrt{2}(\sqrt{1-x}+\sqrt{x})} \geq 0 \Leftrightarrow$$

$$\Leftrightarrow \frac{\left(1-\sqrt{2x}\right)^2(1+\sqrt{2x})\left(2(\sqrt{2x}+\sqrt{x}+\sqrt{2}) - \frac{(1+\sqrt{2x})(1+2x)}{\sqrt{1-x}+\sqrt{x}}\right)}{2\sqrt{2}(\sqrt{1-x}+\sqrt{x})} \geq 0,$$

$$\Leftrightarrow \frac{\left(1-\sqrt{2x}\right)^2(1+\sqrt{2x})\left(2(\sqrt{2x}+\sqrt{x}+\sqrt{2})\sqrt{1-x} + \sqrt{2x} - 1\right)}{2\sqrt{2}(\sqrt{1-x}+\sqrt{x})^2} \geq 0,$$

It is evident that the latter inequality is true for all values of $x \in (0; 1)$, since

$$\sqrt{2}(\sqrt{1-x}+\sqrt{x}) \geq \sqrt{2}\sqrt{1-x+x} = \sqrt{2} > 1.$$

Thus, the inequality (7) is also true for all values of $x$. Hence,

$$\sqrt{1-a} + \sqrt{1-b} + \sqrt{1-c} + \sqrt{1-d} - \sqrt{a} - \sqrt{b} - \sqrt{c} - \sqrt{d} \geq$$

$$\geq -\sqrt{2}(a^2 + b^2 + c^2 + d^2 - 1) = 0.$$

The proof is complete.

**Example 2:** (Russian Olympiad in Mathematics, final round, 2003, $9^{th}$ grade, Problem 6, by S. Berlov). Let $a, b, c$ – positive numbers whose sum equals 1. Proof that:

$$\frac{1}{1-a} + \frac{1}{1-b} + \frac{1}{1-c} \geq \frac{2}{1+a} + \frac{2}{1+b} + \frac{2}{1+c}.$$

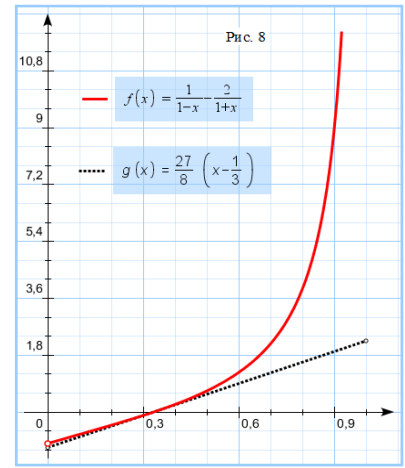

*Solution*: Let $f(x) = \frac{1}{1-x} - \frac{2}{1+x}$, where $x \in (0; 1)$. Let's introduce a tangent line to the graph $f$ at the point $x_0 = \frac{1}{3}$:

$$y = f\left(\frac{1}{3}\right) + f'\left(\frac{1}{3}\right)\left(x - \frac{1}{3}\right) = 0 + \frac{27}{8}\left(x - \frac{1}{3}\right) = \frac{27}{8}\left(x - \frac{1}{3}\right).$$

Since the following (see. pic. 8),

$$\frac{1}{1-x} - \frac{2}{1+x} \geq \frac{27}{8}\left(x - \frac{1}{3}\right)$$

is equivalent to the inequality:

$$(3x - 1)^2(3x + 1) \geq 0,$$

which is true for all values of $x \in (0; 1)$, then the following is also true:

$$\frac{1}{1-a} - \frac{2}{1+a} + \frac{1}{1-b} - \frac{2}{1+b} + \frac{1}{1-c} - \frac{2}{1+c} \geq \frac{27}{8}(a + b + c - 1) = 0.$$

The proof is done.

**Note 5**: From picture 8 it can be seen that the function $f(x) = \frac{1}{1-x} - \frac{2}{1+x}$, where $x \in (0;1)$, is downward convex, however, since the second derivative of the function $f$ is

$$f''(x) = \frac{2}{(1-x)^3} - \frac{4}{(1+x)^3}$$

on the interval (0; 1) takes as positive values as negative values, then the function $f$ is not convex on the interval (0; 1).

**Example 3:** (Saint Petersburg Olympiad, 1988, 9the grade, final round, Problem 3). Let $a, b, c, d$ be positive real numbers. Prove the following:

$$\frac{1}{a}+\frac{1}{b}+\frac{4}{c}+\frac{16}{d} \geq \frac{64}{a+b+c+d}.$$

*Solution*: Since the presented inequality is homogenous, it is sufficient to prove it when $a+b+c+d=8$. In this case, it can be rewritten in the following way:

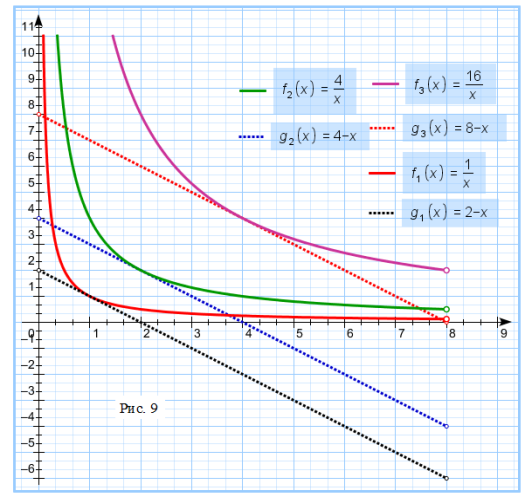

$$\frac{1}{a}+\frac{1}{b}+\frac{4}{c}+\frac{16}{d} \geq 8.$$

Let $f_1(x) = \frac{1}{x}$, $f_2(x) = \frac{4}{x}$, $f_3(x) = \frac{16}{x}$, where $0<x<8$. The points $x_1$, $x_2$, $x_3$ are chosen such that

$$x_1 + x_1 + x_2 + x_3 = 8, f_1'(x_1) = f_2'(x_2) = f_3'(x_3).$$

Рис. 9

It is not hard to note that the following values can be chosen for: $x_1=1$, $x_2=2$, $x_3=4$. Since functions $f_1, f_2, f_3$ are downward convex their graphs lie not below their corresponding tangents (see pic. 9). The following inequalities can be further written:

$$\frac{1}{a} \geq f_1(1) + f_1'(1)(a-1) = 1 - (a-1) = 2 - a,$$

$$\frac{1}{b} \geq f_1(1) + f_1'(1)(b-1) = 1 - (b-1) = 2 - b,$$

$$\frac{4}{c} \geq f_2(2) + f_2'(2)(a-2) = 2 - (c-2) = 4 - c,$$

$$\frac{16}{d} \geq f_3(4) + f_3'(4)(d-4) = 4 - (d-4) = 8 - d.$$

Adding these inequalities, we obtain

$$\frac{1}{a}+\frac{1}{b}+\frac{4}{c}+\frac{16}{d} \geq 16 - (a+b+c+d) = 8.$$

The proof is done.

It is challenging to point out all applications of the presented method in one project paper. We have attempted to show the most interesting and notable cases. Therefore, it is suggested that readers apply this new method to problems from the references cited earlier, as well as the following sources [13], [15], [18], [19].

This paper is prepared under the supervision of Ibragim Ibatulin.

## Conclusion

The following was done in this paper:

- ✓ It was shown that among third degree polynomials there exist infinite number of non-convex functions that comply with Jensen's inequality at a given point;
- ✓ Incorrect proofs were corrected for inequalities from Chetkovski Z. and Suppa E.;
- ✓ A new method to proof inequalities was developed which was called in this paper as "Method of Separating Tangents" and presented at XXI International Mathematics Conference «Mathematics and Education», Cheboksary, Russia. The scale of application of the method can be judged from books of Pham Kim Hung (Secrets in Inequalities) and Chetkovski Z. (Inequalities) where applicable inequalities are seen 20% of the time in which the method can proof about 70% of problems.;
- ✓ The class of functions was significantly expanded where Jensen's inequality is applicable at a given point. Moreover, results were obtained and conditions stated for functions to comply with Jensen's inequality at a given point. Theorems 1 to 4;
- ✓ Theorems 3, 4 and 5 of this paper were highly regarded by Hrabrov A.I. as novelties and their proofs with corresponding examples of inequalities were published on arxiv.org (arxiv.org/abs/1311.4404);
- ✓ Alternative proofs were shown for 10 popular inequality problems presented at various international olympiads and popular mathematics books. Overall, proofs for about 100 inequalities were developed and presented.
- ✓ As of today, the newly developed method is published in 4 articles:
- ✓ Cheboksary, Russia
  "About one method of solving inequalities" by I.Z. Ibatulin, A.N. Lepes,
  // Works of XXI conference «Mathematics. Education» in 2013, (pages. 34–50);
- ✓ Sofia, Bulgaria
  "Application of the method of separating tangents to prove inequalities by Ibatulin I.Z., Lepes A.N.
  // Didactical Modeling: e-journal 2013. URL: http://www.math.bas.bg/omi/DidMod/index.htm. 13 pages (to be published in 2014);
- ✓ Hong Kong, China
  Ibatulin I.Zh., Lepes A.N. Using tangent lines to prove inequalities (part II)
  // Mathematical Excalibur. 2013–2014. V.18 N.4, 6 pages.
- ✓ Moscow, Russia
  Ibatulin I.Zh., Lepes A.N. Using tangent lines to prove inequalities
  // Mathematics in School. 2014. №4, p. 20–24.

The above publications contain alternative proofs for 34 inequalities.